\documentclass[12pt,a4paper]{article}
\usepackage{amsmath}
\usepackage{amsfonts}
\usepackage{amssymb}
\usepackage{xspace}

\newtheorem{theorem}{Theorem}

\renewcommand\th{\theta}
\def\l{\left}  \def\r{\right}
\def\a{\alpha} \def\b{\beta}
\renewcommand\t{\tau}

\newcommand\qed{\hfill$\quad${\rule{3mm}{3mm}}\medskip\\}
\newcommand\LB[1]{\label{#1}}
\newcommand\BE[2]{\begin{#1} #2 \end{#1}}
\newcommand\EQ[2]{\BE{equation}{\LB{#1} #2}}

\newcommand{\ay}{asymptotic\xspace}
\newcommand{\eq}{equation\xspace}

\newcommand\ds{\displaystyle}

\newcommand\f{\varphi}

\newcommand \Lm{\Lambda}
\newcommand\Th{\Theta}

\newcommand\dl{\delta}

\newcommand{\bbN}{\blackboard{N}}

\newcommand{\blackboard}[1]{\mathbb#1}
\newcommand{\iy}{\infty}

\title{ 'Wave type' spectrum of  the  Gurtin-Pipkin \eq of the second order
}

\author{
{S. A. Ivanov\thanks{
State Marine Technical University of St. Petersburg.
{Research of Sergei Ivanov
was supported in part by the Russia Foundation for Basic Research, grant 08-01-00595a.
{\tt sergei.ivanov@pobox.spbu.ru }}}}
}

\begin{document}
\date{}

\maketitle

\begin{abstract}
We study the complex part of the spectrum of the
the Gurtin-Pipkin integral-differential \eq of the second order in time.
We consider the model case when the kernel is a sum of exponentials
$a_k\exp(-b_k)$ with
$a_k=1/k^\a$,  $b_k=k^\b$.
We show that there are two complex
sequences of points of the spectrum asymptotically close to the spectrum points of
the wave \eq.

\end{abstract}

\vskip1cm

\section{\LB{intro_etc} Introduction. Notations. Main theorem}

\subsection{\LB{ss10}  Introduction}

In several fields of physics such as
heat transfer with finite propagation speed \cite{GuPip},
systems with thermal memory \cite{EMV}, viscoelasticity problems \cite{CMD},
and acoustic waves in composite media \cite{Sham},
the following integro-differential \eq arises
\begin{equation}\LB1
\theta_{tt}(x,t)=a \th_{xx}-\int_0^t k(t-s) \th_{xx}\,ds
 \ x\in(0,\pi), \ t>0.
\end{equation}
with the Dirichlet boundary condition and with the initial data
$\th(0,x)=\xi(x)$, $\th_t(0,x)=\eta(x)$.
Regularity of this equation  is studied in \cite{P05} and
in \cite{VW}.

First, apply the Fourier method:
we set $\f_n=\sqrt{\frac2\pi}\sin nx$ and expand the solution and the initial
data  in series in $\f_n$
$$
\th(x,t)=\sum_1^\iy \th_n(t)\f_n(x),\qquad \xi(x)=\sum_1^\iy \xi_n \f_n(x),\
 \eta(x)=\sum_1^\iy \eta_n \f_n(x).
$$
For the components we obtain
\begin{equation}\label{thn}
\ddot \theta_n(t)=-\a n^2\th_n(t)+n^2\int_0^t k(t-s) \theta_n (s)d s, \  t>0,
\ \theta_n(0)=\xi_n,\  \dot\th_n(0)=\eta_n.
\end{equation}

We will denote the Laplace image by the capital characters.
Applying the Laplace Transform to \eqref{thn} we find
$$
z^2\Theta_n(z)-z\xi_n-\eta_n=-a n^2\Th_n(z)+n^2 K(z) \Theta_n(z)
$$
or
$$
\Theta_n(z)=\frac{z\xi_n+\eta_n}{z^2+an^2-n^2K(z)}.
$$

The set $\Lm$ of all zeros of the
denominators
$z^2+an^2-n^2K(z)$, $n=1,2,\dots$
is called the spectrum of the equation \eqref1.

\subsection{\LB{ss11}  The statement and the class of kernels}

In application \cite{Sham} the kernel $k(t)$ is a sum of exponentials
$$
k(t)=\sum_1^\iy a_ke^{-b_kt}.
$$
We consider the model case
$$
a=1, \ a_k=1/k^\a, \ b_k=k^\b, \ \a >0.
$$
The problem is to describe the zeros $z_n$ with $|\Im z_n|\gg1$  of
the functions
\EQ0{
z^2/n^2+1-K(z),\  n\in \bbN,
}
where $K$ is the Laplace transform of $k(t)$ in our case,
$$
K(z)=\sum_1^\iy \frac {a_k}{z+b_k}.
$$
We are interested in the case
\EQ{coeff1}
{
\sum_1^\iy a_k=\iy,\  \sum_1^\iy \frac {a_k}{b_k}<\iy.
}
what gives
$$
\a\le 1, \\a+\b>1.
$$
The conditions \eqref{coeff1} means that $k(t)$ has an integrable singularity at zero:
$k\notin L^\iy(0,\iy)$, $k\in L^1(0,\iy)$.


\subsection{\LB{ss12}  The main result}
Set
$$
r=\frac{\a+\b-1}\b, \ 0<r\le 1.
$$

\begin{theorem}
{\LB {mainth}}
For any $n$ there exist two zeros $z_n^\pm$, $z_n^-=\overline{z_n^+}$ of \eqref0 such that

\par(i) If $r<1$,
$$
z_n^+ =in+c_re^{-i(r+1)\pi/2} n^{1-r} +O(1)+O(n^{1-2r}), \  c_r=\frac{\pi}{\b\sin \pi r}.
$$

(ii)If $r=1$,
$$
z_n^\pm=\pm in-\frac1{2\b}\log n +O(1).
$$
\end{theorem}

\section{\LB{main_steps} The proof of the  main theorem}
We are going to find the zeros 'near' the positive imaginary axis and we set
$$
z_n=in+\t_n n.
$$
Then we obtain the \eq
\EQ{2.1}{
\t_n(\t_n+2i)=K(z_n).
}
Consider the case  $|\t_n|\ll 1$.
\BE{remark}{
If we consider the case where the second factor in\eqref{2.1}
 is small: $|\t_n+2i|\ll 1$ we obtain the complex
conjugated roots. Indeed $\bar z_n=-in+\bar {\t_n} n$ is also the solution to \eqref0.
}
Rewrite \eqref{2.1} as
\EQ{root}{
\t=g_n(\t), \ g_n(\t)=\frac{K(in+\t n)}{\t+2i}.
}
Show that for large $n$ we have a contraction map and we can applied the fixed point theorem.
We need
\BE{lemma}{\LB{zK'}
If for some $\dl>0$
\EQ{dl}{
|\arg z|<\pi-\dl,
}
then for $|z|\to\iy$

(i)
 $r<1$
$$
K(z)=c_r z^{-r}+ O(|z|^{-1}).
$$

(ii)
 $r=1$.
$$
K(z)=\frac1\b\frac{\log(1+z)}z + O(1/|z|).
$$

(iii)
\EQ{est}{
|zK'(z)|\prec \BE{cases}{
1/\rho^r,& \, r<1\\
\frac{\log\rho}{\rho}, & \, r=1.
}
}
}

The proof is in Sec. \ref{auxi}.
Using this lemma, we find the \ay of $\t_n$.
Let $|\t|<1/2$ (and $z=in+n\t$)
$$
|g_n'(\t)|=\l|
\frac{K'(z)n}{\t+2i}-K(z)\frac1{(\t+2i)^2}
\r|
$$
$$
\le
|K'(z)|\frac{n}{|\t+2\i|}+|K(z)|\frac1{|\t+2i|^2}\prec
|K'(z)z|+|K(z)|
$$
$$
\prec \BE{cases}{
1/\rho^r,& \, r<1\\
\frac{\log\rho}{\rho}, & \, r=1.
}
$$
as $|z|\to\iy$.
Therefore we can find $R$ such that $|g'_n(\t)|<\rho<1$ for $|z|>R$ and $|\t|<1/2$. This
implies  that \eqref{root} has a zero $\t_n$ with $|\t_n|<1/2$.

The iterations $\t_n^{(0)}=0,$
$\t_n^{(1)}=g_n(\t_n^0)=g_n(0)=-\frac i2 K(in)$, $\dots,\t_n^{(k+1)}=g_n(\t_n^{(k+1)})$,
\dots,  converge and we may write
$$
\t_n=\t_n^{(1)}+(\t_n^{(2)}-\t_n^{(1)})+\dots+(\t_n^{(k+1)}-\t_n^{(k)})+\dots
=\t_n^{(1)}+O(|\t_n^{(2)}-\t_n^{(1)}|).
$$
Estimate $|\t_n^{(2)}-\t_n^{(1)}|$.
$$
|\t_n^{(2)}-\t_n^{(1)}|=|g_n(\t_n^{(1)})-g_n(\t_n^{(0)})|\le
\max_{|\t|\le |\t_n^{(1)}|}|g_n(\t)||\t_n^{(1)}|
$$
$$
\prec
|zK'(z)|\prec \BE{cases}{
1/\rho^{2r},& \, r<1\\
\frac{\log^2\rho}{\rho^2}, & \, r=1.
}.
$$

Now Theorem \ref{mainth} follows from Lemma \ref{zK'} and from
$$
z_n^+=in+n\t_n=in-\frac{in}{2}K(in)+O(n|\t_n^{(2)}-\t_n^{(1)}|).
$$

\qed

\section{\LB{auxi} Asymptotic of $K(z)$}
\subsection{Integral instead series}
Find the error when we replace the series $K(z)$ by the integral
$$
h(z)=\int_1^\iy \frac{dx}{x^\a(z+x^\b)}.
$$
\BE{lemma}{\LB{est2}
Under \eqref{dl} and $x\ge1$, $b>0$
$$
|z+x^b|^2\asymp |z|^2+x^{2b}.
$$
}
Proof. Set $z=\rho e^{i\f}$, $\rho=|z|$
$$
|z+x^b|^2=(\rho\cos\f+x^b)^2+\rho^2\sin^2\f=\rho^2+2x^b\rho\cos\f+x^{2b}.
$$
Evidently
$$
\rho^2+2x^b\rho\cos\f+x^{2b}
\le
(\rho+x^{b})^2\le 2(\rho^2+x^{2b})
$$
and
$$
\rho^2+2x^b\rho\cos\f+x^{2b}
\ge
\rho^2-2x^b\rho\cos\dl+x^{2b}\ge
(\rho^2+x^{2b})-(\rho^2+x^{2b})\cos\dl
$$
$$
=
(1-\cos\dl)(\rho^2+x^{2b}).
$$
\qed
\BE{lemma}{\LB{series}
\EQ{4_1}{
|K(z)-h(z)|\prec
1/\rho.
}
}
Proof.
Set
$$
w(x)=\frac1{x^\a(z+x^\b)}.
$$
$$
|K(z)-h(z)|=\l|
\sum_1^\iy
\int_k^{k+1}\l[
w(k)-w(x)
\r]dx
\r|
\le
\sum_1^\iy
\int_k^{k+1}\l|w(k)-w(x)\r|dx
$$

We have for $x\in[k,k+1]$
$$
|w(x)-w(k)|\le (x-k)\max|w'(x)|,
$$
$$
|w'(x)|=\l|\frac
{\a x^{\a-1}z+(\a+\b) x^{\a+\b-1}}{x^{2\a}(z+x^\b)^2}
\r|
\asymp\l|\frac{\rho+x^\b}{x^{a+1}(\rho^2+x^{2\b})} \r|,
$$
what gives
$$
|w'|\prec \frac{1 }{k^{\a+1}(\rho+k^\b)}
$$
For $\a>0$ this gives
$$
\rho|K(z)-h(z)|\le\sum_1^\iy \frac1{k^{\a+1}}=Const<\iy.
$$

\qed


\BE{remark}{The restriction $r\le1$ in the theorem connected with the estimate
\eqref{4_1}. If $\rho>1$, then the error $1/\rho$
of replacing the series by the integral is more than the
integral, which h is  $O(1/\rho^r)$.
}

\subsection{Evaluating of the integral}
\BE{lemma}{\LB{asy_int}
(i)
For $r<1$
$$
h(z)=c_r z^{-r}+O(1/|z|)
$$

(ii)
For $r=1$
$$
h(z)=\frac1\b\frac{\log(z+1)}z.
$$
}
Proof.

\textbf{(i)}

First, reduce one parameter. Set $x^\b=t$, $x=t^{1/\b}$, $dx=\frac1\b t^{1/\b-1}$. Then
$$
h(z)=\frac1\b\int_1^\iy \frac{dt}{t^r(z+t)}.
$$
Replace this integral by the integral on $[0,\iy)$ with
the error
$$
\l| \int_0^1 \frac{dt}{t^r(z+t)} \r|\asymp \frac1{|z|}.
$$
Therefore
$$
h(z)=\frac1\b\int_0^\iy \frac{dt}{t^r(z+t)}+O\l(\frac1{|z|}\r)=
\frac1\b J(z)+O\l(\frac1{|z|}\r).
$$

In $J(z)$ we change  the variable $y=t/z$. We obtain the integral by the ray
$R(z)=\{a e^{i\arg z}\,\Big|\, a\ge0\}$.
$$
J(z)=z^{-r}\int_{R(z)}\frac {dy}{y^r(1+y)}.
$$
The integrand is an analytical function in the upper half plane and is $o(1/M)$ on
the circle $|z|=M$. This implies that the integral does not depend on the ray.
Now
$$
J(z)=z^{-r}\int_0^\iy \frac {dt}{t^r(1+t)}= \frac{\pi}{\sin\pi r}z^{-r},
$$
see \cite[Probl. 28.22(7)]{1} or \cite[Probl. 878]{2}

\textbf{(ii)}

This is direct evaluating.
\qed

Proof of Lemma \ref{zK'}.
parts (i) and (ii) follows from Lemmas \ref{zK'} and \ref{series}.
Part (iii).
$$
|zR'(z)|\le\sum_1^\iy\frac{|z|}{k^\a|z+k^\b|^2}\overset{\mbox{Lemma} \, \ref{est2}}
\prec \rho\sum_1^\iy\frac1{k^\a(\rho^2+k^{2\b})}.
$$
Replace the series by integral:
$$
|zK'(z)|\le \rho\frac1{\rho^2+1}+\sum_2^\iy\frac1{k^\a(\rho^2+k^{2\b})}
\prec \frac 1\rho+\int_1^\iy \frac{dx}{x^\a(\rho^2+x^{2\b})}.
$$
Estimate the last integral  setting
$x^\b=\rho t$
$$
J=\int_1^\iy \frac{dx}{x^\a(\rho^2+x^{2\b})}\prec \rho^{1/\b-\a/\b-2}
\int_{1\rho}^\iy \frac{dt}{y^r(1+t^2)}.
$$
This
gives
$$\ds
J\prec\BE{cases}{
1/\rho^r,&\, r<1\\
\frac1\rho \log\rho,&\, r=1.
}
$$

\qed


\textbf{Acknowledgements}

The author is very grateful to Prof. A. E. Eremenko  for his constructive
comments and consultations. The author is  also grateful to
Prof. V. V. Vlasov for the fruitful discussions on the topic.

\end{document}